\documentclass[11pt]{amsart}
\newtheorem{thm}{Theorem}[section]
\newtheorem{lema}{Lemma}[section]
\newtheorem{prop}{Proposition}[section]

\newcommand{\be}{\begin{eqnarray}}
\newcommand{\ee}{\end{eqnarray}}
\newcommand{\ben}{\begin{eqnarray*}}
\newcommand{\een}{\end{eqnarray*}}

\newcommand{\RR}{\mathbb R}
\numberwithin{equation}{section}

\title[Uniqueness of ground states]
      {On the uniqueness of positive solutions of
            a quasilinear equation containing a
                     weighted $p$-Laplacian, the superlinear case }\thanks{This research was supported by
        FONDECYT-1030593.}
\author[M. Garc\'{\i}a-Huidobro
        and
        Duv\'an A. Henao]{Marta Garc\'{\i}a-Huidobro C. and Duv\'an A. Henao M.}

\subjclass{35J70, 35J60} \keywords{positive solutions, uniqueness,
superlinear, separation.}

 \email{mgarcia@mat.puc.cl}
 \email{dahenao@puc.cl}

\begin{document}

\maketitle

  {\small \centerline{ Department of Mathematics }
  \centerline{ Pontificia Universidad Cat\'olica de Chile } \centerline{ Casilla 306, Correo 22,
   Santiago, Chile } }
 \medskip

%%%%%%%%%%%%%%%%%%%%%%%%%%%%%%%%%%%%%%%%%%%%%%%%%%%%%%%%%%%%%%%%%%%%%%%%

\begin{abstract}
We consider the problem of uniqueness of  positive solutions to
$$ -\Delta_p
u=K(|x|)f(u),\quad p>1,  \leqno( P) $$  in $\mathbb R^n$,  where
$\Delta_pu:=\mbox{div}(|\nabla u|^{p-2}\nabla u)$. Here $K$ is a
positive $C^1$ function defined in $\RR^+$ and $f\in C[0,\infty)$
has one zero at $u_0>0$, is non positive and not identically 0 in
$(0,u_0)$, and it is locally lipschitz, positive and satisfies
some superlinear growth assumption in $(u_0,\infty)$.
\end{abstract}

\centerline{\today}

\section{Introduction and main results}

The problem of {\it
uniqueness} of radial ground states of the equation $\Delta_p
u+f(u)=0$ and of various related equations has been studied with great
detail during the last decades, see for example
the works of \cite{cl},  \cite{coff}, \cite{cfe1}, \cite{cfe2},
\cite{fls}, \cite{Kw}, \cite{m}, \cite{ms}, \cite{pel-ser1},
\cite{pel-ser2}, \cite{pu-ser}, and \cite{st} among others.
In \cite{pghms},
 a first step in the study of the
{\it uniqueness} of radial ground states and various qualitative properties of
solutions
of
\begin{eqnarray}\label{pghms-}
-\mbox{div}(A(|x|)|\nabla u|^{p-2}\nabla u)=B(|x|)f(u),\quad x\in\RR^n,\quad n>1,
\end{eqnarray}
where $p>1$ and $A,\ B$ are positive $C^1$ functions defined in $\RR^+$, was  accomplished.
(A nonnegative solution of \eqref{pghms-} which tends to 0 as $|x|\to\infty$
is called a ground state solution).
This was done for the case of a {\em sublinear} $f$ but under considerably mild assumptions on $f$ near the origin.

Very recently, in \cite{cgh}  the case of $p=2$ with a {\em superlinear } $f$ in \eqref{pghms-},
was considered at the cost of imposing
a stronger growth condition on the weight functions and a convexity condition on $f$.
Following the ideas in  \cite{cfe2}, \cite{cgh},  in this paper we will establish
the uniqueness of radially symmetric
ground state solutions of
\begin{eqnarray}\label{eq1}
\begin{gathered}
-\Delta_p u=K(|x|)f(u),\quad x\in \RR^n, \quad n>p>1,\\
u(x)\ge0,\ u\not\equiv0,\quad x\in \RR^n,\quad \lim\limits_{|x|\to\infty}u(x)=0,
\end{gathered}
\end{eqnarray}
as well as the uniqueness of positive solutions to
\begin{eqnarray}\label{eqB}
\begin{gathered}
-\Delta_p u=K(|x|)f(u),\quad x\in B_R(0), \quad n>p>1,\\
u(x)>0,\quad x\in B_R(0),\quad u(x)=0,\quad |x|=R.
\end{gathered}
\end{eqnarray}
We are also able to deal with the corresponding homogeneous
Dirichlet--Neumann free boundary problem
\be\label{2.19}
\begin{gathered}
-\Delta_p u=K(|x|)f(u),\quad x\in B_R(0), \quad n>p>1,\\
u(x)>0,\quad x\in B_R(0),
\quad u=\frac{\partial u}{\partial\nu}=0\quad\mbox{on }\partial B_R(0),
\end{gathered}
\ee
which can be seen as a particular case of either \eqref{eq1} (compact support solution) or \eqref{eqB},
see section \ref{examples}.

We will assume that $K:\RR^+\to\RR^+$ is a positive $C^1$ function satisfying
\begin{enumerate}
\item[$(K_1)$] The function $\displaystyle r\to p+\frac{rK'(r)}{K(r)}$ is strictly
positive and decreasing in $(0,\infty)$,
\end{enumerate}
and that the function $f:[0,\infty)\to\RR$ satisfies
\begin{enumerate}
\item[$(f_1)$]$f(0)=0$, and there exists $u_0>0$ such that $f(u)>0$ for $u>u_0$, and
$f(u)\le 0$, $f(u)\not\equiv 0$, for $u\in(0,u_0)$.
\item[$(f_2)$] $f$ is continuous in $[0,\infty)$ and locally lipschitz continuous in $[u_0,\infty)$,
\end{enumerate}
Under assumption $(f_2)$, $f$ is a.e. differentiable in $[u_0,\infty)$. Denoting by $D$ the subset of $[u_0,\infty)$
where $f'$ exists, we will assume the following superlinearity and convexity conditions:
\begin{enumerate}
\item[$(f_3)$] $(p-1)f(u)\le f'(u)(u-u_0)$, for all $u\in D$,
\item[$(f_4)$] The function $u\to\displaystyle\frac{uf'(u)}{f(u)}$
is decreasing in $D$.
\end{enumerate}
In the case of radially symmetric solutions, problems \eqref{eq1} and \eqref{eqB} take  the forms
\begin{eqnarray}\label{eq2}
\begin{gathered}
-(r^{n-1}\phi_p(u'))'=r^{n-1}K(r)f(u),\quad r>0,\quad n>2,\\
u'(0)=0,\quad \lim\limits_{r\to\infty}u(r)=0,\quad u(r)\ge 0\
\mbox{for }r>0,
\end{gathered}
\end{eqnarray}
and
\begin{eqnarray}\label{eq2B}
\begin{gathered}
-(r^{n-1}\phi_p(u'))'=r^{n-1}K(r)f(u),\quad r>0,\quad n>2,\\
u'(0)=0,\quad u(R)=0,\quad u(r)> 0\ \mbox{for }r\in(0,R),
\end{gathered}
\end{eqnarray}
respectively, where we have denoted $\phi_p(s):=|s|^{p-2}s$, $\phi(0)=0$. Also, by a solution to the equation in
\eqref{eq2}, \eqref{eq2B}, we mean a function $u\in C^1[0,\infty)$ such that $r^{n-1}\phi_p(u')\in C^1[0,\infty)$.
\medskip

We state next our main results.

\begin{thm}\label{main1}
Assume that $f$ satisfies $(f_1)$-$(f_4)$, and that the weight $K$ satisfies $(K_1)$.
Then problem \eqref{eq2} has at most one non trivial solution.
\end{thm}

\begin{thm}\label{main2}
Assume that $f$ satisfies $(f_1)$-$(f_4)$, and that the weight $K$ satisfies $(K_1)$.
Then problem \eqref{eq2B} has at most one non trivial solution.
\end{thm}

Our work will follow the ideas in \cite{cl}, \cite{coff},
\cite{cfe2}, \cite{cgh}, \cite{fls}, \cite{Kw}, \cite{pel-ser1}
and \cite{pel-ser2}. That is, we will consider an initial value
problem associated to the equation in \eqref{eq2} or \eqref{eq2B},
and study the behavior of the solution and its derivative with
respect to the initial value.

The work in \cite{pghms} was done in the case
that the radial version of
\eqref{pghms-}, that is
\begin{eqnarray}\label{pghms-r}
-(a|u'|^{p-2}u')'=b(r)f(u),\quad r>0,\quad n>1,
\end{eqnarray}
(with the obvious notation $a(r)=r^{n-1}A(r)$,  $b(r)=r^{n-1}B(r)$) could be transformed, via a
diffeomorphism $r=r(t)$ of $\RR_0^+$, into the form
\begin{eqnarray}\label{pghms-trans}
-(q(t)|v_t|^{p-2}v_t)_t=q(t)f(v),\quad t>0,
\end{eqnarray}
that is, to a problem having the same positive weight function $q$ at both sides, and that
$$q_t>0,\quad \frac{q_t}{q}\quad\mbox{strictly decreasing for $t>0$ and }\ \lim_{t\to0^+}\frac{tq_t}{q}=N-1\ge 0$$
for some $N\in\mathbb R$.
This requires that the functions $a,\ b$
satisfy the  assumptions

\noindent
$$(b/a)^{1/p}\in L^1(0,1)\setminus L^1(1,\infty),\leqno(W1)$$
$$\quad \Bigl(\frac{1}{p}\frac{a'}{a}+\frac{1}{p'}\frac{b'}{b}\Bigr)
\Bigl(\frac{a}{b}\Bigr)^{1/p}\quad\mbox{is strictly decreasing and positive in $\mathbb R^+$},\leqno(W2)$$
($ p'=p/(p-1)$) and  there is $N\ge 1$ such that
$$\lim_{r\to0^+}\Bigl(\frac{1}{p}\frac{a'}{a}+\frac{1}{p'}\frac{b'}{b}\Bigr)
\Bigl(\frac{a}{b}\Bigr)^{1/p}\int_0^r\Bigl(\frac{b}{a}\Bigr)^{1/p}ds=N-1.\leqno(W3)$$

Clearly, our assumptions $n>p$ and $(K_1)$ imply that  $K^{1/p}\in L^1(0,1)\setminus L^1(1,\infty)$,
$r\to r^pK(r)$ is strictly increasing, implying that $(W1)-(W2)$ are satisfied, and
also that the function $ r\to rK(r)$ is   in $L^1(0,1)$.

We note also that our results will cover the uniqueness of ground
state solutions to \eqref{pghms-r}, when $a,\ b:(0,\infty)\to\RR$
are positive $C^1$ functions satisfying
\begin{enumerate}
\item[$(w_1)$] $a^{1-p'}\in L^1(1,\infty)\setminus L^1(0,1)$, (so
that $h:r\to \int_r^{\infty}a^{1-p'}(t)dt$ is well defined)
\item[$(w_2)$] $\displaystyle
\Bigl(\frac{1}{p}\frac{a'}{a}+\frac{1}{p'}\frac{b'}{b}\Bigr)\frac{h}{|h'|}$
is decreasing in $(0,\infty)$,\\
with
$\displaystyle \lim_{r\to\infty}\Bigl(\frac{1}{p}\frac{a'}{a}+\frac{1}{p'}\frac{b'}{b}\Bigr)\frac{h}{|h'|}\ge p-1$.

\end{enumerate}
Indeed,  we can make the change of variable
$$t:=(h(r))^{-(p-1)/(N-p)},\quad v(t)=u(r),$$
where $N>p$ is arbitrary,
to transform \eqref{pghms-r} into
\begin{eqnarray}\label{trans-eq}
-(t^{N-1}\phi_p(v_t))_t=t^{N-1}\tilde K(t)f(v),\quad t>0,\quad N>p,
\end{eqnarray}
where now
$$\tilde K(t)=(\frac{N-p}{p-1})^pa^{p'-1}(r(t))b(r(t))h^{\frac{N(p-1)}{N-p}+1}(r(t)),$$
thus
$$p+\frac{t\tilde K_t(t)}{\tilde K(t)}=p'\frac{N-p}{p-1}\Bigl(\Bigl(\frac{1}{p}\frac{a'}{a}+
\frac{1}{p'}\frac{b'}{b}\Bigr)\frac{h}{|h'|}
-(p-1)\Bigr).$$
Then, in view of assumption $(w_2)$
we see that $\tilde K$ satisfies $(K_1)$.
\bigskip

Our work is organized as follows. In section \ref{prelim} we
introduce some notation and list some general properties of the
solutions to the initial value problem associated to \eqref{eq2}.
In section \ref{basic-mon}, we define an energy-like functional
and establish a monotonicity type of result, which is the key to
prove our uniqueness theorems. Then in section \ref{sep-proofs},
we establish several monotone separation results and prove
Theorems \ref{main1} and \ref{main2}. In section \ref{examples},
we make a remark concerning the support of a ground state solution
and we give some examples to illustrate the type of weights that
we are considering. Finally, there is an appendix where we prove a
basic result concerning the existence of the derivative of $u$ and
$u'$ with respect to the initial data.

\section{Preliminary results}\label{prelim}

The aim of this section is to establish several properties of the solutions to the initial value problem
\begin{eqnarray}\label{ivp}
\begin{gathered}
(\phi_p(u'))'(r)+\frac{n-1}{r}\phi_p(u'(r))+K(r)f(u)=0
\quad r>0,\quad n>p,\\
u(0)=\alpha\quad u'(0)=0.
\end{gathered}\end{eqnarray}
for $\alpha\in(0,\infty).$
Local existence and uniqueness of  solutions to
\eqref{ivp} under assumptions $(f_1)$, $(f_2)$ and $(K_1)$ (in fact only under the assumption $rK\in L^1(0,1)$
for the weight $K$) can be proved as in the appendix in \cite{pghms}, so
we omit proving it.

We will use the equation in \eqref{ivp}  written in any of the three forms
\begin{eqnarray}
\label{std}(r^{n-1}\phi_p(u'))'&=&-r^{n-1}K(r)f(u) \\
\label{PhiP}\left ( \frac{|u'|^p}{p'} \right )' &=&- u'\left
\{\frac{n-1}{r}\phi_p(u')+K(r)f(u)\right \} \\
\label{u''}(p-1)\frac{u''}{u'}&=&-\frac{n-1}{r}-K(r)\frac{f(u)}{\phi_p(u')}.
\end{eqnarray}

As usual, we denote by $u(r,\alpha)$ a  solution of
\eqref{ivp}. Let us set $F(s)=\int_0^sf(t)dt$ and define the
functional \begin{eqnarray}\label{funct-1}
E(r,\alpha)=\frac{|u'(r,\alpha)|^p}{p'K(r)}+F(u(r,\alpha)). \end{eqnarray}
A
simple calculation yields
\begin{eqnarray}\label{I'}
\frac{\partial}{\partial
r}E(r,\alpha)=-\frac{|u'|^p}{p'rK(r)}\Bigl(\frac{n-1}{p-1}p+\frac{rK'(r)}{K(r)}\Bigr),
\end{eqnarray}
and therefore, as $n>p$, we have that $E$ is decreasing in $r$.
Also, $\lim\limits_{r\to 0^+}\frac{|u'|^p}{K(r)}$ exists and we
may apply L'H\^opital's rule to compute it:
\begin{eqnarray*}
\lim_{r\to 0^+}\frac{|u'|^{p-1}}{K^{1/p'}(r)}=\lim_{r\to
0^+}\frac{r^{n-1}|u'|^{p-1}}{r^{(n-1)}K^{1/p'}(r)} =\lim_{r\to
0^+}\frac{p'rK^{1/p}(r)f(u(r))}{\frac{n-1}{p-1}p+\frac{rK'(r)}{K(r)}}=0.
\end{eqnarray*}
It follows then that if $\alpha\le u_0$, then $E(r,\alpha)<0$ for all $r>0$, and thus there cannot exist $0<R\le\infty$
such
that  $u(R,\alpha)=0$. Since we are interested in solutions to \eqref{eq2} or \eqref{eq2B},
we will assume that $\alpha>u_0$.
\medskip

It can be seen that for $\alpha\in(u_0,\infty)$, one has $u(r,\alpha)>0$ and $u'(r,\alpha)<0$ for
$r$ small enough, so  we can define
$$R(\alpha):=\sup\{r>0\ |\ u(t,\alpha)>0\mbox{ and }u'(t,\alpha)<0\ \mbox{ for all }t\in(0,r)\}.$$

Hence, for a given $\alpha\in(u_0,\infty)$, there is a unique solution $u(r,\alpha)$ defined in
$[0,R(\alpha))$, and this function and its derivative are continuous functions of $\alpha$. Also, $u(r,\alpha)$
is invertible in $[0,R(\alpha))$ and we denote its inverse by $t(s,\alpha)$.

Some immediate properties of the solution of \eqref{ivp} are
the following:

\noindent (i) $\lim\limits_{r\to R(\alpha)}u(r,\alpha)$ exists. (We denote by $u(R(\alpha),\alpha)$ this limit).
\medskip

\noindent (ii)  $u(R(\alpha), \alpha)\leq u_0$.
\medskip

\noindent (iii) If $R(\alpha)<\infty$, or if $R(\alpha)=\infty$ and
$u(R(\alpha),\alpha)<u_0$, then $\displaystyle \lim_{r\to
R(\alpha)} u'(r,\alpha)$ exists and it is finite. (We denote by $u'(R(\alpha),\alpha)$ this limit).  In the second
case it is also satisfied that $\displaystyle \lim_{r\to\infty}
r|u'(r,\alpha)|=0$.

\medskip

The first property is clear, so let us first prove  (ii). Assume by contradiction that
$u(R(\alpha),\alpha)=\gamma>u_0$. Then
\[r^{n-1}|u'(r,\alpha)|^{p-1}=\int_0^r t^{n-1}
K(t)f(u(t,\alpha))dt \geq Cf(\gamma) r^{n} K(r) \]
for some
positive constant $C$ and $r$ sufficiently large, where we have
used that
\[\lim_{r\to\infty} \frac{\displaystyle \int_0^r t^{n-1} K(t) dt}{r^n K(r)} \lim_{r\to\infty}
\frac{r^{n-1}K(r)}{nr^{n-1}K(r)+r^nK'(r)}=\lim_{r\to\infty}
\frac{1}{n+r\frac{K'(r)}{K(r)}}>0\] since $n>p$ and
$p+r\frac{K'(r)}{K(r)}$ is strictly positive and decreasing in
$(0,\infty)$. Hence, as $r^pK(r)$ increases,
$r^{p-1}|u'(r,\alpha)|^{p-1}\geq C_0>0$ for some positive constant
$C_0$, contradicting the integrability of $u'$ near infinity.

In order to establish (iii), we observe that if $0\leq u(r,\alpha) < u_0$ for $0<r_1\leq r < R$ then from
(\ref{std}) we have that
\[-(r^{n-1}|u'|^{p-2}u')'=r^{n-1}K(r)f(u)\leq 0\quad \textrm{for
$r_1\leq r<R$,}\] implying that $r^{n-1}|u'|^{p-1}$ decreases for
$r_1\leq r<R$ and hence $\displaystyle \lim_{r\to R}
r^{n-1}|u'|^{p-1}$ exists. We may call it $\lambda$, $\lambda\geq
0$. As $n>p$, when $R=\infty$ we have
\[\lim_{r\to\infty} r|u'| =\lim_{r\to\infty}\frac{\Big
(r^{n-1}|u'|^{p-1}\Big)^{p'-1}}{r^\frac{n-p}{p-1}}=0,\]
implying
also that $u'\to0$. It only remains to examine
the case  $u(r,\alpha)>u_0$ for all
$r\in(0,R(\alpha))$. From (\ref{std}) we have that if
$r<R(\alpha)$ then
\[r^{n-1}|u'(r,\alpha)|^{p-1} = \int_0^r t^{n-1}K(t)f(u(t,\alpha))dt\] This is
not possible if $R(\alpha)<\infty$ since the right hand side of this
equation would be strictly positive and $u'(R(\alpha),\alpha)=0$,
by definition of $R(\alpha)$. We have then that $R(\alpha)=\infty$
if $u(R(\alpha),\alpha)\geq u_0$.

\bigskip

Following \cite{pel-ser1}, \cite{pel-ser2} we define now the sets
\begin{eqnarray*}
{\mathcal N}&=&\{\alpha>u_0\ :\ u(R(\alpha),\alpha)=0\quad\mbox{and}\quad u'(R(\alpha),\alpha)<0\}\\
{\mathcal G}&=&\{\alpha>u_0\ :\ u(R(\alpha),\alpha)=0\quad\mbox{and}\quad u'(R(\alpha),\alpha)=0\}\\
{\mathcal P}&=&\{\alpha>u_0\ :\ u(R(\alpha),\alpha)>0\}.
\end{eqnarray*}

If $\alpha\in\mathcal N$, we will refer to $u(\cdot,\alpha)$ as a {\em crossing } solution, and if
$\alpha\in\mathcal G$, we will refer to $u(\cdot,\alpha)$ as a {\em ground state} solution.

Minor changes in the proof of \cite[Lemma 2.1]{cgh} yield:

\begin{lema}\label{open-sets}
${\mathcal N}$ and ${\mathcal P}$ are open sets.
\end{lema}

\bigskip

Further properties of the solutions are the following:

If $\alpha\in\mathcal N$ then $R(\alpha)<\infty$ and
$E(R(\alpha),\alpha)>0$ so \[E(r,\alpha)=E(R(\alpha),\alpha)+
\int_r^{R(\alpha)}
\frac{|u'(t,\alpha)|^p}{p'tK(t)}\left\{\frac{n-1}{p-1}p+t\frac{K'(t)}{K(t)}\right\}dt>0\quad
\textrm{for $r\in(0,R(\alpha))$}.\]

Let now $\alpha\in \mathcal G$. Since $\displaystyle \lim_{r\to
R(\alpha)} r|u'(r)|=0$ and $r^pK(r)$ increases, we conclude that
$\displaystyle \lim_{r\to R(\alpha)} E(r,\alpha)=0$  and thus
\[E(r,\alpha)=\int_r^{R(\alpha)}
\frac{|u'(t,\alpha)|^p}{p'tK(t)}\left\{\frac{n-1}{p-1}p+t\frac{K'(t)}{K(t)}\right\}dt>0\quad
\textrm{for $r\in(0,R(\alpha))$}.\]

Finally, if $\alpha\in \mathcal P$ we have that
$E(R(\alpha),\alpha)=F(u(R(\alpha),\alpha))<0$ since
$\displaystyle \lim_{r\to R(\alpha)}\frac{|u'(r,\alpha)|^p}{K(r)}$
exists and equals zero. Indeed, this limit exists because
$E(r,\alpha)$ has a limit when $r\to R(\alpha)$. This limit is clearly
equal to zero if $R(\alpha)<\infty$. In the case that $R(\alpha)=\infty$ and using that
\[\lim_{r\to\infty} \frac{|u'(r)|^p}{K(r)}=\lim_{r\to\infty}
\frac{r^p|u'|^p}{r^pK(r)},\] we conclude, as $r^pK(r)$ is
increasing and converges to a positive limit or to infinity when
$r\to\infty$, that it must be equal to $0$. (Otherwise,
$\liminf\limits_{r\to\infty}r|u'(r)|>0$,  contradicting the
integrability of $u'$ near infinity).

 Let $r_0(\alpha)$ denote the (unique) value of $r$ such that
 $u(r,\alpha)=u_0$. Clearly, $r_0(\alpha)$ is finite for $\alpha\in{\mathcal G}\cup{\mathcal N}$,
 and it could be that $r_0(\alpha)=\infty$ for some values of $\alpha\in{\mathcal P}$.
 The following lemma will be proved in the appendix:

 \begin{lema}\label{varphi-der} Under assumptions $(K_1)$,
 $(f_1)$ and $(f_2)$, and for every $\alpha\in(u_0,\infty)$, the
 functions $u(r,\alpha)$ and $r^{n-1}\phi_p(u'(r,\alpha))$ (actually $K^{-1/p'}(r)\phi_p(u'(r,\alpha))$) are of class
 $C^1$ in
 $${\mathcal O}=\{(r,\alpha)\ :\
 \alpha\in(u_0,\infty)\quad\mbox{and}\quad
 r\in[0,t(u_0,\alpha)) \}.$$
For $r\in (0,t(u_0,\alpha))$, we set
$$\varphi(r,\alpha)=\frac{\partial u}{\partial\alpha}(r,\alpha)$$
Then, $\varphi$ satisfies the linear
differential equation
\begin{eqnarray}\label{varphi-eq}
\begin{gathered}
(p-1)(r^{n-1}|u'(r,\alpha)|^{p-2}\varphi'(r,\alpha))'+r^{n-1}K(r)f'(u)\varphi(r,\alpha)=0,\quad
r>0,
\quad n>p,\\
\varphi(0,\alpha)=1,\quad \varphi'(0,\alpha)=0,\quad
\lim\limits_{r\to0^+}r^{n-1}|u'(r,\alpha)|^{p-2}\varphi'(r,\alpha)=0,
\end{gathered}\end{eqnarray}
 and the relation $\varphi'(r,\alpha)=\frac{\partial u'}{\partial\alpha}(r,\alpha)$, $r\ge 0$.
 Also, in case that $r_0(\alpha)<\infty$, both $\varphi(r,\alpha)$ and
 $r^{n-1}|u'(r,\alpha)|^{p-2}\varphi'(r,\alpha)$ can be
 extended continuously to
$$\underline{{\mathcal O}}=\{(r,\alpha)\ :\
 \alpha\in(u_0,\infty)\quad\mbox{and}\quad
 r\in[0,t(u_0,\alpha)] \},$$
 and  the extension of $\varphi'(r,\alpha)$ is the left
 derivative with respect to $r$ of $\varphi(r,\alpha)$.
\end{lema}

\section{Basic monotonicity results}\label{basic-mon}
In order to prove
our results, we divide our analysis  into two parts. First, it is shown that if
$\alpha_1\in\mathcal{G}\cup\mathcal N$ and $\alpha_2>u_0$ are such that the corresponding solutions $u_1$, $u_2$ satisfy
\begin{itemize}
\item $t_2(u_0)\leq t_1(u_0)$ \item
$t_2(u_0)|u_2'\big(t_2(u_0)\big)|>t_1(u_0)|u_1'\big(t_1(u_0))|$,
\end{itemize}
where $t_1$ and $t_2$ denote the inverses of $u_1$ and $u_2$,
respectively, then $R(\alpha_2)<R(\alpha_1)$ and $R(\alpha_2)u_2'(R(\alpha_2))<R(\alpha_1)u_1'(R(\alpha_1))$, where
the last expression is to be understood as a limit when $R(\alpha_1)=\infty$.
This implies in particular, that $\alpha_2\in\mathcal N$, and $u_2'(R(\alpha_2))<u_1'(R(\alpha_1))$.

Analogously, if $\alpha_1\in\mathcal G$ and
\begin{itemize}
\item $t_2(u_0)\geq t_1(u_0)$ \item
$t_2(u_0)|u_2'\big(t_2(u_0)\big)|<t_1(u_0)|u_1'\big(t_1(u_0)\big)|$,
\end{itemize}
then $u_2(R(\alpha_2))>0$ and thus $\alpha_2\in\mathcal P$.

This directs our attention towards the function
$r_0(\alpha)=t(u_0, \alpha)$,  and to the function $ru'$
when $u$ reaches $u_0$. The second part of this work consists in proving that $r_0(\alpha)$ is decreasing and that
$r_0(\alpha)|u'(r_0(\alpha),\alpha)|$ increases with $\alpha$ near
$\bar \alpha$, when $\bar \alpha$ is the initial value of a ground
state or a crossing solution. It shall be seen that
\[\frac{\partial
r_0(\alpha)}{\partial \alpha}
=\frac{\varphi(r_0(\alpha),\alpha)}{|u'(r_0(\alpha),\alpha)|}\]
and
\[r_0(\alpha)^{\frac{n-p}{p-1}-1} \frac{\partial}{\partial \alpha}
r_0(\alpha)u'(r_0(\alpha), \alpha) = \frac{\partial}{\partial r}
\left \{r^{\frac{n-p}{p-1}}\varphi(r,\alpha) \right
\}_{r=r_0(\alpha)}\] Therefore, it is natural to analyze carefully
the zeros of the function $\varphi(r,\alpha)$ for $\alpha$ near
$\bar \alpha$.

These two results combined together lead us to a monotone
separation type of result that will  finally yield the desired
uniqueness theorems.

As mentioned in \cite{cfe2}, the function $\varphi(r,\alpha)$ and
the idea of studying its zeros first appeared in the work of
\textsc{Coffman} \cite{coff}. As we are studying the function
$ru'$,  the function $r|u'(r,\alpha)|/u(r,\alpha)$, which was
proven by \textsc{Kwong} \cite{Kw} to be increasing in $r$ for
$r\in(0,r_0(\bar\alpha))$, will play an important role.

Finally we mention that since all the differential equations that we are given are of the form
$$(r^{n-1}|u'|^{p-2}v')'=\xi(r)$$
for an appropriate function $\xi$, and $v$ any among $u,\ u',\
\varphi$, the following
 identity
$$(r^{n-1}|u'|^{p-2}(v'w-vw'))'=(r^{n-1}|u'|^{p-2}v')'w-(r^{n-1}|u'|^{p-2}w')'v$$
will prove to be crucial as it will allow us to obtain information about $v'w-vw'$, which
will be especially relevant when any of
$v,\ v',\ w,$ or $w'$ is equal to $0$.

\subsection{Behavior of the solutions below $u_0$}

\begin{prop} \label{bajou_0}
Let $\alpha_1\in\mathcal{G}\cup\mathcal{N}$,
and let $\alpha_2>u_0$. Let $t_1$ and $t_2$ denote the
 inverses of $u(r,\alpha_1)$, and $u(r,\alpha_2)$ respectively. If
\[t_2(u_0)\leq t_1(u_0)<\infty \quad\textrm{and}\quad
t_2(u_0)|u_2'(t_2(u_0))|>t_1(u_0)|u_1'(t_1(u_0))|\] then
$u_2(R(\alpha_2))=0$, and
\[t_2(s)<t_1(s)\quad \textrm{and}\quad
t_2(s)|u_2'(t_2(s))|>t_1(s)|u_1'(t_1(s))|\quad \textrm{for $s\in
[0, u_0)$}.\] Conversely, if
\[t_1(u_0)\leq t_2(u_0)<\infty \quad\textrm{and}\quad
t_2(u_0)|u_2'(t_2(u_0))|<t_1(u_0)|u_1'(t_1(u_0))|\] then
\[t_2(s)>t_1(s)\quad \textrm{and}\quad
t_2(s)|u_2'(t_2(s))|<t_1(s)|u_1'(t_1(s))|\quad \textrm{for $s\in
[u_2(R(\alpha_2)), u_0)$}.\]
\end{prop}
\noindent{\bf Remark.} We note that when $s=0$ in the first
statement of the proposition, the conclusion should be read as
$R(\alpha_2)<\infty$ and $R(\alpha_2)|u_2'(R(\alpha_2))|>0$.
\begin{proof}
Let us define
\[\overline F(s)= \int_0^s
f(\xi)t_1^p(\xi)K\big(t_1(\xi)\big)d\xi\] and
\[I(s,\alpha)=t^p(s,\alpha)\frac{|u'(t(s,\alpha),\alpha)|^p}{p'}+\overline F(s) \]

We observe first that $\overline F(s)$ is well defined: indeed,
since $\displaystyle \lim_{r\to R(\alpha_1)} ru_1'(r)= 0$, it
follows that $r^{p-1}|u_1'(r)|^{p-1}$ is bounded and thus
$r^{p-1}|u_1'(r)|^{p-1}|u'(r)|\in L^1(0,R(\alpha_1))$. Using now
that from equation (\ref{PhiP}) we have that
\begin{eqnarray*}
\int_0^s f(\xi)t_1^p(\xi)K\big(t_1(\xi)\big)d\xi &=&
\int_{t_1(s)}^{R(\alpha_1)} r^pK(r)f(u_1(r))|u_1'(r)|dr \\
&=& \left.r^p\frac{|u'|^p}{p^\prime}\right|_{t_1(s)}^{R(\alpha_1)}
+ (n-p)\int_{t_1(s)}^{R(\alpha_1)}r^{p-1}|u_1'(r)|^{p-1}|u'(r)|dr,
\end{eqnarray*}
our claim  follows.

 Using (\ref{PhiP}) again we have that
\begin{eqnarray*}
\frac{\partial}{\partial s}I(s,\alpha) &=&
(p-1)t'(s,\alpha)t^{p-1}(s,\alpha)|u'(t(s,\alpha),\alpha)|^p\qquad\qquad\\
&& +t^p(s,\alpha)
\frac{\partial}{\partial r} \left
\{\frac{|u'(r,\alpha)|^p}{p'}\right \}_{r=t(s,\alpha)}
t'(s,\alpha) +f(s)t_1^p(s)K\big (t_1(s)\big )\\
&=&(p-1)t^{p-1}(s,\alpha)\phi_p(u')-t^p(s,\alpha) \left
\{\frac{n-1}{t(s,\alpha)} \phi_p(u')+f(s)K\big(t(s,\alpha)\big)
\right \} \\
&&\qquad\qquad+f(s)t_1^p(s)K\big(t_1(s)\big)\\
&=&-(n-p)t^{p-1}(s,\alpha)
\phi_p(u')+f(s)\{t_1^p(s)K\big(t_1(s)\big)-t^p(s,\alpha)K\big(t(s,\alpha)\big)\}
\end{eqnarray*}
Evaluating at $\alpha=\alpha_1$ we have:
\[\frac{\partial}{\partial s}
I(s,\alpha_1)=(n-p)t_1^{p-1}(s)|u_1'\big(t_1(s)\big)|^{p-1}>0\]
Besides,
\[I(0,\alpha_1)=R(\alpha_1)^p\frac{|u_1'(R(\alpha_1))|^p}{p'}+0\ge0\]
(if $R(\alpha_1)=\infty$ this continues to be true), and hence
\begin{equation}
\label{I>0}I(s,\alpha_1)>0\quad\forall\,s\in(0,\alpha_1)\end{equation}

Let us assume, for contradiction, that there exists $s_1\in
[0,u_0)$ such that $t_2(s)<t_1(s)$ for all $s \in (s_1, u_0)$ and
$t_2(s_1)=t_1(s_1)<\infty$ (in case that $t_2(u_0)=t_1(u_0)$,
since
\[t_2(u_0)|u_2'(t_2(u_0))|>t_1(u_0)|u_1'(t_1(u_0))|\] then
$t_2(s)$ would be strictly less than $t_1(s)$ for $s$ in a
neighborhood below $u_0$). When $s\to s_1^+$ we have
\[\frac{t_2(s)-t_2(s_1)}{s-s_1} < \frac{t_1(s)-t_1(s_1)}{s-s_1}\]
so $t_2'(s_1)\leq t_1'(s_1)$ and
\[t_2(s_1)|u_2'(t_2(s_1))|=-\frac{t_2(s_1)}{t_2'(s_1)} \leq
-\frac{t_1(s_1)}{t_1'(s_1)}=t_1(s_1)|u_1'(t_1(s_1))|\] (the case
$t_1'(s_1)=-\infty$ is also being considered here). Thus the
existence of $\bar s\in [0,u_0)$ such that
\begin{equation}
\label{relImp}
\begin{array}{cl}
I(s,\alpha_2) > I(s,\alpha_1)>0 & \textrm{for}\ s\in (\bar s,u_0]\\
I(\bar s,\alpha_2)=I(\bar s, \alpha_1)\geq 0 \\
t_2(s)\leq t_1(s) &\textrm{for}\ s\in [\bar s,u_0]
\end{array}
\end{equation}
follows. If instead of the previous case we have that
$t_2(s)<t_1(s)$ for all $s\in (0,u_0)$, $u_1$ and $u_2$ are
defined in $(0,\infty)$, and both $u_1(r)$ and $u_2(r)$ tend to
zero when $r\to\infty$, then \[\lim_{s\to0^+}
I(s,\alpha)=\lim_{r\to\infty} r|u'(r,\alpha)|=0\] for both
$\alpha=\alpha_1$ and $\alpha=\alpha_2$, so we may choose $\bar s$
as the value of $s$ where
\[I(\xi,\alpha_2)>I(\xi, \alpha_1)\quad\textrm{for all $\xi\in
[s,u_0]$}\] ceases to be true. This can also be done if $\alpha_2$
is not in $\mathcal G$ nor in $\mathcal N$ and $t_2(s)<t_1(s)$ for
$s\in(u_2(R(\alpha_2),u_0)$, as in this case we have that
$$t_2(u_2(R(\alpha_2)))|u_2'(R(\alpha_2))|=0<t_1(u_2(R(\alpha_2)))|u_1'(R(\alpha_2))|.$$
Finally, we observe that if $\alpha_2\in\mathcal G\cup\mathcal N$,
$t_2(s)<t_1(s)$ for all $s\in (0,u_0)$ and condition
\[t_2(s)|u_2'(t_2(s))|>t_1(s)|u_1'(t_1(s))|\quad \textrm{for $s\in
[0, u_0)$}\] is not satisfied, then again it can be found $\bar s$
such that (\ref{relImp}) holds.

From relation (\ref{relImp}) we see that we can well define, for
$s\in (\bar s,u_0]$, $W(s,\alpha)=I(s,\alpha)^{\frac{1}{p}}$ for
$\alpha=\alpha_1$ and $\alpha=\alpha_2$, that
$W(s,\alpha_2)-W(s,\alpha_1)$ is positive at $s=u_0$ and equals
zero at $s=\bar s$. Nevertheless this is not possible because
$W(s,\alpha_2)-W(s,\alpha_1)$  decreases for $s\in[\bar s,u_0]$.
Indeed,
\begin{eqnarray*}
\lefteqn{\frac{\partial}{\partial s}
\{W(s,\alpha_2)-W(s,\alpha_1)\} = \frac{1}{p} \left
\{I(s,\alpha_2)^{-\frac{1}{p'}}\frac{\partial}{\partial
s}I(s,\alpha_2) - I(s,\alpha_1)^{-\frac{1}{p'}}
\frac{\partial}{\partial s} I(s, \alpha_1) \right \}}\\
&=& \frac{n-p}{p} \left \{
\frac{t_2^{p-1}(s)|u_2'\big(t_2(s)\big)|^{p-1}}{\displaystyle
\left (t_2^p(s) \frac{|u_2'\big(t_2(s)\big)|^p}{p'}+\overline F(s)
\right
)^{\frac{1}{p'}}}-\frac{t_1^{p-1}(s)|u_1'\big(t_1(s)\big)|^{p-1}}{\displaystyle
\left (t_1^p(s) \frac{|u_1'\big(t_1(s)\big)|^p}{p'}+\overline F(s)
\right )^{\frac{1}{p'}}}\right \}\\
&&+\frac{1}{p}\frac{f(s)}{I(s,\alpha_2)^\frac{1}{p'}} \left (
t_1^p(s)K\big(t_1(s)\big)-t_2^p(s)K\big(t_2(s)\big) \right ).
\end{eqnarray*}

Concerning the first term, using the identity $p-1=\frac{p}{p'}$
we have \ben \frac{[t_2|u_2'|]^{p-1}}{\displaystyle \left
(\frac{1}{p'}[t_2|u_2'|]^p+\overline F\right )^\frac{1}{p'}} -
\frac{[t_1|u_1'|]^{p-1}}{\displaystyle \left
(\frac{1}{p'}[t_1|u_1'|]^p+\overline F\right )^\frac{1}{p'}}\qquad\qquad\qquad\\
\qquad\qquad\qquad=\left (\frac{1}{p'}+\frac{\overline
F}{t_2^p|u_2'|^p}\right )^{-\frac{1}{p'}}-\left
(\frac{1}{p'}+\frac{\overline F}{t_1^p|u_1'|^p}\right
)^{-\frac{1}{p'}}\leq 0 \een since $t_2|u_2'|>t_1|u_1'|$ and
$\overline F(s)\leq 0$ for $s<u_0$.

With respect to the second term, as $r^pK(r)$ is increasing and
since $t_2(s)\leq t_1(s)$ we will have that $\displaystyle f(s)
\left ( r^pK(r) \big |_{r=t_2(s)}^{t_1(s)}\right )\leq 0$ because
$f(s)\leq 0$ for $s\leq u_0$. This completes the proof of the first
part of the proposition.

In order to proof the second part of the proposition assume for a
contradiction that there exists $\bar s\in [u_2(R(\alpha_2)),u_0)$
such that
\begin{equation}
\label{relImp2}
\begin{array}{cl}
I(s,\alpha_2) < I(s,\alpha_1) & \textrm{for}\ s\in (\bar s,u_0]\\
I(\bar s,\alpha_2)=I(\bar s, \alpha_1) \\
t_2(s)\geq t_1(s) &\textrm{for}\ s\in [\bar s,u_0]
\end{array}
\end{equation}

Since \[\frac{\partial}{\partial s}I(s,\alpha_2)
=(n-p)t_2^{p-1}(s)|u_2'(s)|^{p-1}
+f(s)\{t_1^p(s)K\big(t_1(s)\big)-t_2^p(s)K\big(t_2(s)\big)\}>0\]
for $s\in (\bar s, u_0)$ and $I(\bar s, \alpha_2)=I(\bar s,
\alpha_1)\geq 0$ then $W(s,\alpha_2)$ is well defined for $s$ in
that interval. From relation (\ref{relImp2}) we obtain that
 $W(s,\alpha_2)-W(s,\alpha_1)$ would be negative at $s=u_0$
and equal to zero at $s=\bar s$, and this is not possible since
$W(s,\alpha_2)-W(s,\alpha_1)$ is increasing in $[\bar s,u_0]$.
Indeed,
\begin{eqnarray*}
\lefteqn{\frac{\partial}{\partial s}
\{W(s,\alpha_2)-W(s,\alpha_1)\} = \frac{1}{p} \left
\{I(s,\alpha_2)^{-\frac{1}{p'}}\frac{\partial}{\partial
s}I(s,\alpha_2) - I(s,\alpha_1)^{-\frac{1}{p'}}
\frac{\partial}{\partial s} I(s, \alpha_1) \right \}}\\
&=& \frac{n-p}{p} \left \{
\frac{t_2^{p-1}(s)|u_2'\big(t_2(s)\big)|^{p-1}}{\displaystyle
\left (t_2^p(s) \frac{|u_2'\big(t_2(s)\big)|^p}{p'}+\overline F(s)
\right
)^{\frac{1}{p'}}}-\frac{t_1^{p-1}(s)|u_1'\big(t_1(s)\big)|^{p-1}}{\displaystyle
\left (t_1^p(s) \frac{|u_1'\big(t_1(s)\big)|^p}{p'}+\overline F(s)
\right )^{\frac{1}{p'}}}\right \}\\
&&+\quad\frac{1}{p}\frac{f(s)}{I(s,\alpha_2)^\frac{1}{p'}} \left (
t_1^p(s)K\big(t_1(s)\big)-t_2^p(s)K\big(t_2(s)\big) \right ).
\end{eqnarray*}
The first term will be positive since $\overline F(s)<0$ and
$t_2|u_2'|<t_1|u_1'|$, and the second term will be positive since
$t_2(s)\geq t_1(s)$ and $f(s)\leq 0$. This completes the proof.
\end{proof}

\subsection{Behavior of the solutions above $u_0$}

We  begin the second part of this work by noting that if $u$
satisfies (\ref{std}) and we set \[v(r):=ru'(r)+cu(r),\] where $c$
is an arbitrary constant, then using (\ref{u''}) we have
\begin{equation}
\label{v'}
v'(r)=\left(c-\frac{n-p}{p-1}\right)u'(r)-\frac{1}{p-1}rK(r)\frac{f(u)}{|u'(r)|^{p-2}}
\end{equation}

Using now (\ref{std}) we have
\begin{eqnarray}
\nonumber (p-1)(r^{n-1}|u'|^{p-2}v')'&=& \nonumber -
(c(p-1)-(n-p))r^{n-1}K(r)f(u)-(r^nK(r)f(u))'\\
&=& -r^{n-1}K(r)f(u) \left
\{c(p-1)+\left(p+r\frac{K'(r)}{K(r)}\right)+ru'(r)\frac{f'(u)}{f(u)}\right\}\label{v}
\end{eqnarray}

Next we prove the following lemma. In what follows we denote
$r_0(\alpha)$ for $r(u_0,\alpha)$.

\begin{lema} \label{elLema} Let
$\bar\alpha\in\mathcal{G}\cup\mathcal{N}$.
 Then
$ru'(r)/u(r)$ is strictly decreasing in $(0,\bar r_0)$, where
$u(r):=u(r,\bar\alpha)$ and $\bar r_0:=r_0(\bar\alpha)$.
\end{lema}

\begin{proof}
As in \cite{cfe2} and \cite{cgh}, we set $w(r)=ru'(r)$. Then
\[u^2\left (\frac{ru'}{u}\right)'(r)=uw'(r)-wu'(r)\] and it
suffices to show that $wu'-uw'$ is positive in $(0,\bar r_0)$. Let
$r\in(0,\bar r_0)$. It is easily seen that \ben
(p-1)r^{n-1}|u'|^{p-2}
(wu'(r)-uw'(r))\qquad\qquad\qquad\\
=\int_0^r\{(p-1)(t^{n-1}|u'|^{p-2}u')'w-(p-1)(t^{n-1}|u'|^{p-2}w')'u\}dt
\een
and since $w(r)$ is a particular case of the function $v(r)$
defined above (where $c=0$),
\[(p-1)(r^{n-1}|u'|^{p-2}w')'=-r^{n-1}K(r)f(u) \left
\{\left(p+r\frac{K'(r)}{K(r)}\right)+w(r)\frac{f'(u)}{f(u)}\right\}.\]
Therefore,
\begin{eqnarray*}
(p-1)r^{n-1}|u'|^{p-2}
(wu'(r)-uw'(r))
=\int_0^rt^{n-1}K(t)w(t)(u(t)f'(u)-(p-1)f(u))dt\\
\qquad\qquad\qquad+\int_0^rt^{n-1}K(t)f(u)u(t)\left(p+\frac{tK'(t)}{K(t)}\right)dt\\
=\int_0^r t^nK(t)((uf(u))'-pf(u))dt
+\int_0^rt^{n-1}K(t)f(u)u(t)\left(p+t\frac{K'(t)}{K(t)}\right)dt
\end{eqnarray*}

Let us define \[F_0(s)=\int_{u_0}^s f(\xi)d\xi,\] then we have
\begin{eqnarray*}\lefteqn{\int_0^rt^nK(t)((uf(u))'-pf(u))dt}\\
&=&r^nK(r)(uf(u)-pF_0(u))-
\int_0^rt^{n-1}K(t)uf(u)\left(1-p\frac{F_0(u)}{uf(u)}\right)\left(n+t\frac{K'(t)}{K(t)}\right)dt
\end{eqnarray*}
 and hence
\begin{eqnarray}
(p-1)r^{n-1}|u'|^{p-2} (wu'(r)-uw'(r))\label{kw}
 = r^nK(r)(uf(u)-pF_0(u))\qquad\qquad\nonumber\\
\qquad\qquad +\int_0^rt^{n-1}K(t)uf(u)\left\{
p\frac{F_0(u)}{uf(u)}\left(n+t\frac{K'(t)}{K(t)}\right)-(n-p)\right
\}dt
\end{eqnarray}

Clearly, the hypothesis $f'(u)(u-u_0)\geq(p-1)f(u)$ for $u>u_0$
implies that
\begin{equation}\label{F_0}pF_0(u)\leq f(u)(u-u_0)\ \textrm{if $u\geq u_0$}\,,\end{equation}
so it remains only to show that
\begin{equation}\label{aux2}\int_0^rt^{n-1}K(t)uf(u)\left\{
p\frac{F_0(u)}{uf(u)}\left(n+t\frac{K'(t)}{K(t)}\right)-(n-p)\right
\}dt>0\quad \textrm{for $r\in(0,\bar r_0)$}\end{equation}

Set
\[G(s)=p\left(n+t(s,\bar\alpha)\frac{K'(t(s,\bar\alpha))}{K(t(s,\bar\alpha))}\right)\frac{F_0(s)}{sf(s)}-(n-p)\]

We see first that evaluating at $r=\bar r_0$ equation (\ref{kw})
becomes
\begin{equation} \label{aux1}
(p-1)\bar r_0^{n-1}|u'(\bar r_0)|^{p-2}(wu'(\bar r_0)-uw'(\bar
r_0))= \int_0^{\bar r_0}
t^{n-1}K(t)u(t)f(u(t))G(u(t))dt\end{equation}

On the other hand, using (\ref{v'}) we obtain that $w'(\bar
r_0)=-\frac{n-p}{p-1}u'(\bar r_0)$, that is, \[wu'(\bar
r_0)-uw'(\bar r_0)=-|u'(\bar
r_0)|\left(\frac{n-p}{p-1}u_0+ru'(\bar r_0) \right )\] Using
(\ref{v'}) once again we obtain
\[\frac{\partial}{\partial r}\left \{\frac{n-p}{p-1}u(r)+ru'(r)\right\}=-
\frac{rK(r)}{p-1}\frac{f(u)}{|u'|^{p-2}}\geq 0\] for $r\in (\bar
r_0, R(\bar\alpha))$, and since $\bar \alpha\in \mathcal G\cup
\mathcal N$ we must have \[\lim_{r\to R(\bar \alpha)}
\frac{n-p}{p-1}u(r)+ru'(r)\leq 0,\] implying
that\[\left(\frac{n-p}{p-1}u_0+ru'(\bar r_0) \right )<0\] so
$wu'(\bar r_0)-uw'(\bar r_0)>0$ and from (\ref{aux1}) we see that
$G$ cannot be everywhere negative in $(u_0, \bar \alpha)$.
Relation (\ref{F_0}) shows us that
\[0\leq \lim_{s\to u_0^+}\frac{pF_0(s)}{f(s)} \leq \lim_{s\to
u_0^+}(s-u_0)=0\] and since $n>p$ and $n+t(s,\bar \alpha)
\frac{K'(t(s,\bar \alpha))}{K(t(s,\bar \alpha))}$ is strictly
positive we know that $G(u_0)=-(n-p)<0$. Finally, we observe that
\[\left(\frac{F_0(s)}{sf(s)}\right)'=(s^2f(s))^{-1} \left
\{sf(s)-F_0(s)\left (1+\frac{sf'(s)}{f(s)}\right)\right\}\] By
hypothesis, $\frac{sf'(s)}{f(s)}$ decreases when $s$ increases, so
$sf(s)-F_0(s)\left (1+\frac{sf'(s)}{f(s)}\right)$ is increasing in
$s$ and equals zero when $s=u_0$, implying that
$\frac{F_0(s)}{sf(s)}$ is increasing in $s$. $n+t(s,\bar \alpha)
\frac{K'(t(s,\bar \alpha))}{K(t(s,\bar \alpha))}$ is increasing
too, since $n>p$, so $G(s)$ is increasing in $s$. We conclude the
existence of $s_1\in (u_0,\bar\alpha)$ such that $G(s_1)=0$,
$G(s)<0$ for $s\in (u_0, s_1)$ and $G(s)>0$ for $s\in(s_1, \bar
\alpha)$. Now, if $u(r)\geq s_1$ then (\ref{aux2}) follows. If
$u(r)<s_1$ then
\begin{eqnarray*}
\lefteqn{\int_0^rt^{n-1}K(t)uf(u)\left\{
p\frac{F_0(u)}{uf(u)}\left(n+t\frac{K'(t)}{K(t)}\right)-(n-p)\right
\}dt}\\
&>& \int_0^{\bar r_0} t^{n-1}K(t)uf(u)\left\{
p\frac{F_0(u)}{uf(u)}\left(n+t\frac{K'(t)}{K(t)}\right)-(n-p)\right
\}dt \\
&=& (p-1)\bar r_0^{n-1}|u'(\bar r_0)|^{p-2}(wu'(\bar r_0)-uw'(\bar
r_0))>0,
\end{eqnarray*}
and thus the lemma follows.
\end{proof}

\section{Monotone separation results and proof of the main
theorems}\label{sep-proofs}

We begin this section with our first monotone separation result.

\begin{thm}\label{sobreu_0}
If $\bar \alpha \in \mathcal G\cup \mathcal N$ then $r_0(\alpha)$
is decreasing and $r_0(\alpha)u'(r_0(\alpha), \alpha)$ is strictly
decreasing in a neighborhood of $\bar\alpha$.
\end{thm}

\begin{proof}
From the relation $u(r_0(\alpha),\alpha)=u_0$ we obtain
\begin{equation}\label{r_0'}\frac{\partial r_0}{\partial \alpha}= -
\frac{\varphi(r_0(\alpha),\alpha)}{u'(r_0(\alpha),\alpha)}.\end{equation}
On the other hand, using equation (\ref{u''}) (which is simplified
because it is evaluated at $r=r_0(\alpha)$ so that the term
containing $f(u)$ disappears) and computing directly it follows
that
\begin{equation}\label{energia}
\frac{\partial}{\partial \alpha}\{
r_0(\alpha)u'(r_0(\alpha),\alpha)\}=\frac{n-p}{p-1}\varphi(r_0(\alpha),\alpha)
+r_0(\alpha)\varphi'(r_0(\alpha),\alpha).
\end{equation}

For $r\in(0,r_0(\alpha))$ we have
\begin{eqnarray*}
\lefteqn{(p-1)(r^{n-1}|u'|^{p-2}\varphi')'(u-u_0)=\{(p-1)r^{n-1}|u'|^{p-2}\varphi'(u-u_0)\}'}
\\
&&\qquad\qquad\qquad\qquad\qquad\qquad\qquad -(p-1)\varphi'r^{n-1}\phi_p(u')\\
&=&(p-1)\{r^{n-1}|u'|^{p-2}\varphi'(u-u_0)-\varphi
r^{n-1}\phi_p(u')\}'+(p-1)\varphi(r^{n-1} \phi_p(u'))'\\
&=&(p-1)\{r^{n-1}|u'|^{p-2}[\varphi'(u-u_0)-\varphi
u']\}'-(p-1)r^{n-1}Kf\varphi,
\end{eqnarray*}
and therefore we have from equation (\ref{varphi-eq}) that
\begin{equation}\label{aux3}
(p-1)\{r^{n-1}|u'|^{p-2}[\varphi'(u-u_0)-\varphi
u']\}'+r^{n-1}K\varphi[f'(u)(u-u_0)-(p-1)f(u)]=0.
\end{equation}
In particular
\begin{equation}\label{aux4}
(p-1)r_0^{n-1}|u'(r_0)|^{p-1}
\varphi(r_0,\alpha)+\int_0^{r_0}r^{n-1}K\varphi[f'(u)(u-u_0)-(p-1)f(u)]dr=0,
\end{equation}
and this is valid for every $\alpha$ (not only for $\bar \alpha$).
For simplicity, we will write $r_0=r_0(\alpha)$, $\bar
r_0=r_0(\bar \alpha)$ and $\varphi (r)=\varphi(r,\bar \alpha)$. We
know, from (\ref{varphi-eq}) that $\varphi(0)=1$. Assume that
$\varphi(r)>0$ for $r\in(0,\bar r_0)$. Since, by hypothesis
$(f_3)$, $f'(u)(u-u_0)-(p-1)f(u)\geq 0$ for $u\in D$, \eqref{aux4}
yields  $\varphi(\bar r_0)=0$ (implying immediately that
$\varphi'(\bar r_0)<0$) and  $f'(s)(s-u_0)-(p-1)f(s) \equiv 0$ for
$s\in (u_0, \bar \alpha)$. Replacing in (\ref{energia}) we have,
in particular, that $\frac{\partial}{\partial
\alpha}\{r_0(\alpha)u'(r_0(\alpha),\alpha)\}<0$ at $\alpha=\bar
\alpha$. From (\ref{aux4}) we deduce that $\varphi(r_0, \alpha)=
0$ when $\alpha$ is close to $\bar \alpha$ but less than
$\bar\alpha$, and $\varphi(r_0, \alpha)\le 0$ if
$\alpha>\bar\alpha$ because $u(r,\alpha)\in (u_0,\bar \alpha)$
except  for $r\in(0,t(\bar\alpha,\alpha))$, where
$\varphi(r,\alpha)$ is positive as $\varphi(0,\alpha)=1$ for all
$\alpha$ and $\varphi$ is continuous. Hence we conclude that
$\frac{\partial r_0}{\partial \alpha}\leq
0$ in a neighborhood of $\bar \alpha$.\\

The case when $\varphi$ has a first zero $r_1$ in $(0,\bar r_0)$ is
more difficult and must be treated differently. The idea is to
show that $r_1$ is the only zero of $\varphi$ in $(0,\bar r_0]$,
implying that $\varphi(\bar r_0)<0$ (which immediately yields
$\frac{\partial r_0}{\partial \alpha} (\bar \alpha)<0$), and that
also $\frac{n-p}{p-1} \varphi(\bar r_0, \bar \alpha) + \bar r_0
\varphi'(\bar r_0, \bar \alpha)<0$, yielding again
$\frac{\partial}{\partial
\alpha}\{r_0(\alpha)u'(r_0(\alpha),\alpha)\}<0$. Here we will make
use of our Lemma \ref{elLema}.\\

Let us set as before $v(r)=ru'(r)+cu(r)$, where $u(r)$ denotes
$u(r,\bar \alpha)$. From equation (\ref{v}) we see that $v$
satisfies
\begin{equation}\label{v2}(p-1)(r^{n-1}|u'|^{p-2}v')'+r^{n-1}K(r)f'(u)v =r^{n-1}\Omega(r)\end{equation}

where \[\Omega(r)
=-K(r)f(u)\left\{c(p-1)+\left(p+r\frac{K'(r)}{K(r)}\right)-c\frac{uf'(u)}{f(u)}\right
\}.\]

\noindent $\Omega(r)$ can also be written as
\[\Omega(r)=-K(r)f(u)\left
\{\left(p+r\frac{K'(r)}{K(r)}\right)-c\frac{f'(u)(u-u_0)-(p-1)f(u)}{f(u)}-cu_0\frac{f'(u)}{f(u)}\right\}\]

\noindent $\Omega(r_1)$ would be negative if $c=0$; By hypothesis,
\[c\frac{f'(u(r_1))(u(r_1)-u_0)-(p-1)f(u(r_1))}{f(u(r_1))}\ge 0\] if $c>0$
because $r_1\in(0,\bar r_0)$; and since
$u_0\frac{f'(u(r_1))}{f(u(r_1))}$ is positive, the expression
\[p+r_1\frac{K'(r_1)}{K(r_1)}-cu_0\frac{f'(u(r_1))}{f(u(r_1))}\]
can be made negative if we choose $c$ large enough. Therefore
there exists $c>0$ such that $\Omega(r_1)=0$. Since, by the
hypotheses imposed on $K$ and on $f$
\[c(p-1)+\left(p+r\frac{K'(r)}{K(r)}\right)-c\frac{u(r)f'(u(r))}{f(u(r))}\]
is decreasing in $r$, we have that $\Omega(r)$ is negative in
$(0,r_1)$ and positive in $(r_1, \bar r_0)$.\\

Using relations (\ref{v2}) and (\ref{varphi-eq}) we see that
\begin{eqnarray*}
\lefteqn{(p-1)r^{n-1}|u'|^{p-2}(\varphi v'(r)-v\varphi'(r))}\\
&=&(p-1)\int_0^r\{(t^{n-1}|u'|^{p-2}v')'\varphi-(t^{n-1}|u'|^{p-2}\varphi')'v\}dt\\
&=& \int_0^r\{-t^{n-1}K(t)f'(u)v(t)+t^{n-1}\Omega(t)\}\varphi(t)dt +
\int_0^r t^{n-1}K(t)f'(u)\varphi(t)v(t)dt\\
&=&\int_0^r t^{n-1}\Omega(t)\varphi(t) dt
\end{eqnarray*}

We conclude therefore that
\begin{equation}\label{clave}
\varphi(r)v'(r)-v(r)\varphi'(r)\leq 0
\end{equation}
for all $r$ in $(0,r_1]$, and also for $r$ in $(r_1, \bar r_0]$ as
long as $\varphi(r)$ remains negative in $(r_1, r)$. In particular
\eqref{clave} is true for $r=r_1$, so we have that $v(r_1)$ is
necessarily non positive. Now we can show that $\varphi(r)$ has no
zeroes in $(r_1,\bar r_0]$, because if otherwise there would be a
first zero $r_2$ after $r_1$, for which we should have, according
to (\ref{clave}), that $v(r_2)\ge0$. But since $ru'(r)/u(r)$ is
strictly decreasing in $(0,\bar r_0)$,
\[v(r)=u(r) \left (
\frac{ru'(r)}{u(r)}+c\right)<u(r)\left (
\frac{r_1u'(r_1)}{u(r_1)}+c\right)=\frac{u(r)}{u(r_1)}v(r_1)<0\
\textrm{for $r>r_1$}.\]

We can finally prove our theorem. We have that $\varphi(\bar r_0)<0$
so we already have proven the first part of the theorem. If we had
also that $\varphi'(\bar r_0)<0$ then we would be ready, so lets
assume $\varphi'(\bar r_0)>0$. This can only happen if
$c<\frac{n-p}{p-1}$. Indeed, if $c\geq \frac{n-p}{p-1}$, then, as
can be seen in equation (\ref{v'}) we would have
\[v'(r)=\left(c-\frac{n-p}{p-1}\right)u'(r)-
\frac{1}{p-1}rK(r)\frac{f(u)}{|u'(r)|^{p-2}}\leq
-\frac{1}{p-1}rK(r)\frac{f(u)}{|u'(r)|^{p-2}}\leq 0\] for all
$r\in (0,\bar r_0)$. Since $v(\bar r_0)<0$, if $\varphi'(\bar
r_0)>0$ relation (\ref{clave}) would be contradicted.\\

Hence we shall assume that $c<\frac{n-p}{p-1}$ and that
$\varphi'(\bar r_0)>0$. Since (\ref{clave}) holds for all $r\in
(0,\bar r_0)$, by using (\ref{v'}) we have that \[\left
\{\varphi(r)|u'(r)|\left(\frac{n-p}{p-1}-c\right)-\varphi'(r)v(r)\right\}-\frac{rK(r)}{p-1}
\frac{f(u)\varphi(r)}{|u'|^{p-2}}\leq 0\] hence letting $r\to \bar
r_0$ we obtain that \[\varphi(\bar r_0)|u'(\bar r_0)|\left
(\frac{n-p}{p-1}-c\right)-\varphi'(\bar r_0)v(\bar r_0) \leq 0. \]

On the other hand,
\begin{eqnarray*}
\frac{n-p}{p-1}v(\bar r_0)&=&c \big(\bar r_0u'(\bar
r_0)+\frac{n-p}{p-1}u(\bar r_0)
\big)-\left(c-\frac{n-p}{p-1}\right)\bar r_0u'(\bar r_0) \\
&<& -\left(c-\frac{n-p}{p-1}\right) \bar r_0u'(\bar r_0)
\end{eqnarray*}
because $\bar r_0u'(\bar r_0)+\frac{n-p}{p-1}u(\bar r_0)$ is
negative when $\bar \alpha\in \mathcal G\cup \mathcal N$, as was
proved in Lemma \ref{elLema}. Therefore we have
\begin{eqnarray*}
\varphi(\bar r_0)|u'(\bar r_0)|\left (\frac{n-p}{p-1}-c\right)&\leq &
\varphi'(\bar r_0)v(\bar r_0)\\
&<&-\frac{p-1}{n-p}\varphi'(\bar
r_0)\left(c-\frac{n-p}{p-1}\right) \bar r_0u'(\bar r_0)
\end{eqnarray*}
and this implies, dividing by $\displaystyle \left
(\frac{n-p}{p-1}-c\right)|u'(\bar r_0)|$, that
\[\frac{n-p}{p-1}\varphi(\bar r_0)+\bar r_0 \varphi'(\bar r_0) < 0,\] as we
required.
\end{proof}

 The above result is sufficient in order to prove Theorem \ref{main1}, but thanks to Proposition
 \ref{bajou_0}, it can be extended to the
 following stronger monotone
separation theorem.

\begin{thm} \label{monotonicity1}
If $\bar \alpha\in \mathcal G\cup \mathcal N $ then there exists
$\delta >0$ such that for
$\alpha\in(\bar\alpha,\bar\alpha+\delta)$,
$u(R(\alpha),\alpha)=0$,
$$t(s,\alpha)<t(s,\bar\alpha)\quad\mbox{and}\quad
t(s,\alpha)|u'(t(s,\alpha),\alpha)|>
t(s,\bar\alpha)|u'(t(s,\bar\alpha),\bar\alpha)|$$ for every fixed
$s\in [0, u_0)$.
 If instead, $\alpha\in(\bar\alpha-\delta,\bar\alpha)$,
then
$$t(s,\alpha)>t(s,\bar\alpha)\quad\mbox{and}\quad
t(s,\alpha)|u'(t(s,\alpha),\alpha)|<
t(s,\bar\alpha)|u'(t(s,\bar\alpha),\bar\alpha)|$$ for every fixed
$s\in [u(R(\alpha),\alpha), u_0)$.
\end{thm}

\begin{proof}
We know from Theorem \ref{sobreu_0} that there exists $\delta>0$
such that $r_0(\alpha)$ is decreasing and
$r_0(\alpha)|u'(r_0(\alpha),\alpha)|$ is strictly increasing in
$(\bar \alpha -\delta, \bar \alpha+\delta)$. Therefore for
$\alpha\in (\bar \alpha,\bar \alpha+\delta)$ we have that
\[t(u_0,\alpha)\leq t(u_0,\bar \alpha)<\infty\quad\textrm{and}\quad
 t(u_0,\alpha)|u'(t(u_0,\alpha),\alpha)|>t(u_0,\bar
 \alpha)|u'(t(u_0,\bar\alpha),\bar \alpha)|,\] and if $\alpha\in(\bar\alpha-\delta,\bar\alpha)$, then
\[t(u_0,\alpha)\geq t(u_0,\bar \alpha)\quad\textrm{and}\quad
 t(u_0,\alpha)|u'(t(u_0,\alpha),\alpha)|<t(u_0,\bar
 \alpha)|u'(t(u_0,\bar\alpha),\bar \alpha)|.\]
The result follows now from Proposition  \ref{bajou_0}.
\end{proof}

We deduce now the following separation theorem that will allow us
to prove the uniqueness of radial ground states solutions to
(\ref{eq1}).

\begin{thm} \label{separation}
If $\bar \alpha \in \mathcal G$ then there exists $\delta >0$ such
that \[(\bar \alpha, \bar \alpha+\delta)\subset \mathcal
N\quad\textrm{and}\quad (\bar \alpha-\delta, \bar \alpha)\subset
\mathcal P.\]
\end{thm}

\begin{proof}
We have from the previous theorem that, for $\alpha\in (\bar
\alpha, \bar \alpha+\delta)$, $u(R(\alpha),\alpha)=0$,
\[R(\alpha)=t(0,\alpha)<t(0,\bar\alpha)\le\infty\] and
\[R(\alpha)|u'(R(\alpha),\alpha)|> \lim_{r\to R(\bar \alpha)}
r|u'(r,\bar\alpha)|=0,\] so $\alpha$ must be in $\mathcal N$.

We prove that $(\bar \alpha-\delta, \bar \alpha)\subset \mathcal
P$ by contradiction: let $\alpha\in (\bar \alpha-\delta, \bar
\alpha)$ and assume that $u(R(\alpha),\alpha)=0$. From the
previous theorem we obtain then that
\[R(\alpha)|u'(R(\alpha),\alpha)|<\lim_{r\to R(\bar \alpha)}
r|u'(r,\bar\alpha)|=0,\] which clearly cannot be true.
\end{proof}

\begin{proof}[$\underline{\textit{Proof of Theorem \ref{main1}}}$.]\ Assume by contradiction that
there are two different ground state solutions, an let
$\alpha_1<\alpha_2$ be their initial values. Let us set \[\bar
\alpha =\sup \{\alpha\in (\alpha_1,\alpha_2):\, (\alpha_1, \alpha)
\subset \mathcal N\},\] which is well defined because of the last
theorem applied to $\alpha_1$. $\bar \alpha$ cannot be in
$\mathcal N$ nor in $\mathcal P$ since both sets are open, and
hence $\bar \alpha \in \mathcal G$, but this contradicts the fact
that, by the last theorem applied to $\bar \alpha$, there exists a
neighborhood below $\bar \alpha$ entirely contained in $\mathcal
P$. This completes the proof.
\end{proof}
In order to prove Theorem \ref{main2}, we establish next the
following monotonicity result.

\begin{thm}\label{monotd}
If $\bar \alpha \in \mathcal G\cup\mathcal N$ then $\alpha\in
\mathcal N$ for $\alpha>\bar \alpha$ and $R(\alpha)=t(0,\alpha)$
is strictly decreasing in $[\bar \alpha, \infty)$.
\end{thm}

\begin{proof}
By Theorem \ref{monotonicity1} there exists $\delta>0$ such that
$R(\alpha)$ is strictly decreasing and
$R(\alpha)|u'(R(\alpha),\alpha)|$ is strictly increasing for
$\alpha\in[\bar \alpha,\bar \alpha+\delta)$. This implies, as in
the proof of Theorem \ref{separation}, that $(\bar \alpha, \bar
\alpha+\delta)\subseteq \mathcal N$. Let us set \[\alpha_1=\sup
\{\alpha>\bar \alpha:\, (\bar \alpha, \alpha) \subset \mathcal
N\}.\] If $\alpha_1$ were finite, $\alpha_1$ cannot be in
$\mathcal N$ nor in $\mathcal P$ because of the openness of both
sets, and cannot be in $\mathcal G$ because of Theorem
\ref{separation}, we conclude that $\alpha_1$ must be infinite.
The monotonicity of $R(\alpha)$ for $\alpha\in[\bar\alpha,\infty)$
follows  using Proposition \ref{bajou_0} and the compactness of
$[\alpha_1,\alpha_2]$, $\alpha_1\ge\bar\alpha$, completing the
proof.
\end{proof}

\begin{proof}[$\underline{\textit{Proof of Theorem \ref{main2}.}}$]
It follows directly from Theorem \ref{monotd}.
\end{proof}

Finally, we note that if $\alpha\in\mathcal N$, we can write
$$E(R(\alpha),\alpha)=\frac{t(0,\alpha)^p|u'(t(0,\alpha))|^p}{p't(0,\alpha)^pK(t(0,\alpha))},$$
which is strictly increasing in $\alpha$ by Theorems
\ref{monotonicity1} and \ref{monotd}. This is consistent with the
fact that $E(R(\alpha),\alpha)$ is negative for $\alpha\in
\mathcal P$, equals zero when $\alpha\in \mathcal G$ and is
positive for $\alpha\in \mathcal N$, considering that the initial
values in $\mathcal P$ are below and initial values in $\mathcal
N$ are above the initial value of the ground state solution, in
case it exists. These properties suggest that it can be considered
that $E(R(\alpha),\alpha)$ quantifies how crossing a solution is,
and how far  it is from the unique ground state solution.

\section{Concluding remarks and examples}\label{examples}

\bigskip

We end this article with some remarks and examples that illustrate our result.

First of all, we mention that Theorem 5.7 and Corollary 5.8 in \cite{pghms} concerning the support of ground states,
continue to hold in this case. That is, if $u(\cdot,\bar\alpha)$ is a ground state, then
$R(\bar\alpha)<\infty$ if and only if
$$\int_0\frac{du}{\sqrt{|F(u)|}}<\infty.$$
Indeed, assumption $(K_1)$ implies that all assumptions needed in \cite{pghms} for the weights
to prove this result are satisfied,
and since
these results are independent of the superlinear growth of $f$, they follow.

We end this section with some examples. A typical example of an equation of the form \eqref{eq2},
is the Matukuma equation, namely
\begin{eqnarray}\label{1es}
\Delta u + \frac{f(u)}{1+r^\sigma}=0,
 \qquad n>2,\qquad\sigma > 0,
\end{eqnarray}

Another example is the equation
\begin{eqnarray}\label{2es}
\Delta u+\frac{r^{\sigma}}{(1+r^{2})^{\sigma/2}}
\cdot\frac{f(u)}{r^{2}}=0,
 \qquad n>2,\qquad\sigma>0.
\end{eqnarray}
Equation \eqref{2es} was
first introduced in \cite{bfh},  as a
model of stellar structure.
As a main example, we consider the following equation
which includes as
special cases both \eqref{1es} and \eqref{2es}:
\begin{eqnarray}\label{**}
\begin{gathered}
\mbox{div}(|x|^k |Du|^{p-2}Du) + |x|^\ell
\left(\frac{|x|^s}{1+|x|^s}\right)^
{\sigma/s}f(u) = 0,\quad \mbox{in }\mathbb R^N,\ N>1, \\
k \in \RR,\qquad \ell \in \RR, \qquad s > 0,\qquad \sigma > 0.
\end{gathered}
\end{eqnarray}
Here
$$
a(r)=r^{N+k-1},
\qquad b(r)=r^{N+\ell-1}\left(\frac{r^s}{1+r^s}\right)^{\sigma/s}.
$$
We claim that $(w1)$--$(w2)$ are satisfied if
\begin{eqnarray}\label{conds}
N+k>p\qquad\mbox{and}\qquad \ell \ge k - p.
\end{eqnarray}
Indeed, under the first in \eqref{conds} we have that
$a^{1-p'}(s)=s^{\frac{N+k-1}{1-p}}$ belongs to
$L^1(1,\infty)\setminus L^1(0,1)$ and thus $(w_1)$ is satisfied.
In this case
$$h(r)=\int_r^\infty s^{\frac{N+k-1}{1-p}}ds=\frac{p-1}{N+k-p}r^{\frac{N+k-p}{1-p}}$$
 Also, from the
second in \eqref{conds}, we obtain that $N-1+ \frac{k}{p} +
\frac{\ell}{p'}\ge N+k-p$, and thus using the first in
\eqref{conds} we obtain that
$$\left[\dfrac{1}{p}\dfrac{a'}{a}+\dfrac{1}{p'} \dfrac{b'}{b}\right]
\frac{h}{|h'|}=  \left(N - 1 + \frac{k}{p} + \frac{\ell}{p'} +
\frac{\sigma}{p'}\cdot
\frac{1}{1 + r^s}\right)\frac{p-1}{N+k-p}.
$$
is decreasing as $\sigma>0, \ s>0$, and tends to
$$\left(N - 1
+ \frac{k}{p} + \frac{\ell}{p'} \right)\frac{p-1}{N+k-p}$$ as
$r\to\infty$. Since this quantity is greater than or equal to
$p-1$, we see that $(w_2)$ is satisfied.

Finally we mention that our assumptions also cover weights of the
form
$$K(r)=r^{\theta}\mbox{exp}\Bigl({\frac{-|\log(r)|^2}{2}}\Bigr)\ \mbox{near $r=0$},\qquad K(r)=r^{\theta}\log^{\alpha}(1+r),$$
with $\theta\ge -p$ and $\alpha>0$. Indeed, $(K_1)$ is clearly verified since
$$p+\frac{rK'(r)}{K(r)}=p+\theta+|\log(r)|\quad\mbox{near $r=0$},$$
$$p+\frac{rK'(r)}{K(r)}=p+\theta+\alpha\frac{r}{(1+r)\log(1+r)}$$
respectively,
are decreasing functions. It should be noted that in the first case,
$$\lim\limits_{r\to0^+}\frac{rK'(r)}{K(r)}=\infty.$$
On the other hand, for the second example
$$\lim_{r\to\infty}\frac{rK'(r)}{K(r)}=\theta,$$
hence $(K_1)$ is satisfied for $\theta\ge -p$ and $\alpha>0$.

As for the nonlinearities $f$ that we can cover, we mention the canonical example
$$f(u)=u^{q_1}-u^{q_2}, \quad\mbox{with }0<q_2<p-1\le q_1.$$

\section{Appendix}

{\bf Existence of $\frac{d}{d\alpha}u$ and
$\frac{d}{d\alpha}r^{n-1}\phi_p(u')$.} This section is devoted to
the proof of Lemma \ref{varphi-der}. The proof of this lemma is
very similar to the corresponding one given in \cite{cfe2}, but
since it involves delicate computations due to the degeneracy of
the operator $\Delta_p$ when $p\not=2$, (because $u'(0)=0$), we
give it in detail. Let $r\in[0,r^*]$, with $r^*<r_0(\alpha)$. In
order to prove Lemma \ref{varphi-der}, we use the change of
variables introduced in \cite{pghms}: Set
$$t(r):=\int_0^r K^{1/p}(s)ds,\quad v(t):=u(r).$$
By assumption $(K_1)$, $K^{1/p}\in L^1(0,1)\setminus
L^1(1,\infty)$, thus $t(0)=0$, $t(\infty)=\infty$, and
$r^{n-1}K^{1/p'}(r)\phi_p(v_t(t))=r^{n-1}\phi_p(u'(r))$, where
from this point on we use the subscript $_t$ to denote the
derivatives with respect to $t$. Hence, setting
$q(t):=r^{n-1}(t)K^{1/p'}(r(t))$, we see that $v$ satisfies
\eqref{pghms-trans}, $v(0)=\alpha$, and $v_t(0)=0$. Indeed, from
assumption $(K_1)$, \be\label{a00}
\frac{q_t(t)}{q(t)}=\frac{1}{p'}\Bigl(\frac{n-1}{p-1}p+\frac{r(t)K'(r(t))}{K(r(t))}\Bigr)\frac{1}{r(t)K^{1/p}(r(t))}>0,
\ee hence as $n>p$ we have that $q_t>0$ in $(0,\infty)$. Thus,
from the equation \eqref{pghms-trans}, using first that
$\lim\limits_{t\to0^+}q(t)\phi_p(v_t(t))=\lim\limits_{r\to0^+}r^{n-1}\phi_p(u'(r))=0$,
we find that for any $t>0$, \be\label{a0}
\phi_p(|v_t(t)|)=\frac{1}{q(t)}\int_0^tq(s)f(v(s))ds\le
\frac{f(\alpha)}{q(t)}\int_0^tq(s)ds \ee and thus
$\phi_p(|v_t(t)|)\le f(\alpha)t$, implying that $v_t(0)=0$.
Consequently, $v$ satisfies the initial value problem
\begin{eqnarray}\label{pghms-trans2}
\begin{gathered}
-(q(t)|v_t|^{p-2}v_t)_t=q(t)f(v),\quad t>0,\\
v(0)=\alpha,\quad v_t(0)=0.
\end{gathered}
\end{eqnarray}
Also, for any $0\le t\le t^*=t(r^*)$, \be\label{a1}
|v_t(t)|=\Bigl(\frac{1}{q(t)}\int_0^tq(s)f(v(s))ds\Bigr)^{1/(p-1)}
\ge C^*\Bigl(\frac{1}{q(t)}\int_0^tq(s)ds\Bigr)^{1/(p-1)}, \ee
where $C^*=(f(v(t^*)))^{1/(p-1)}>0$.

Let $\{h_n\}\subset\mathbb R$ be any sequence such that
$\lim\limits_{n\to\infty}h_n=0$, and let us set \be
\Theta_n(t,\alpha)&:=&\frac{\phi_p(v_t(t,\alpha+h_n))-\phi_p(v_t(t,\alpha))}{h_n},\label{def-thetan}\\
\psi_n(t,\alpha)&:=&\frac{v_t(t,\alpha+h_n)-v_t(t,\alpha)}{h_n},\label{def-psin}\\
\varphi_n(t,\alpha)&:=&\frac{v(t,\alpha+h_n)-v(t,\alpha)}{h_n}.\label{def-varphin}
\ee Clearly, \be\label{a2}
-\Theta_n(t,\alpha)=\frac{1}{q(t)}\int_0^tq(s)\frac{f(v(s,\alpha+h_n))-f(v(s,\alpha))}
{v(s,\alpha+h_n)-v(s,\alpha)}\varphi_n(s,\alpha)ds \ee and
\be\label{a22} \varphi_n(t,\alpha)=1+\int_0^t\psi_n(s,\alpha)ds.
\ee From the mean value theorem, we see that there is
$\lambda\in[0,1]$ such that
$$
|\Theta_n(t,\alpha)|=(p-1)(\lambda|v_t(t,\alpha+h_n)|+(1-\lambda)|v_t(t,\alpha)|)^{p-2}|\psi_n(t,\alpha)|,$$
hence, using \eqref{a0} and \eqref{a1}, we conclude that for all $t\in(0,t^*]$,
\be\label{a3}
\qquad C_1\Bigl(\frac{Q(t)}{q(t)}\Bigr)^{(p-2)/(p-1)}|\psi_n(t,\alpha)|\le |\Theta_n(t,\alpha)| \le
C_2\Bigl(\frac{Q(t)}{q(t)}\Bigr)^{(p-2)/(p-1)}|\psi_n(t,\alpha)|
\ee
for some positive constants $C_1,\ C_2$, where we have set $Q(t):=\int_0^tq(s)ds$. Using now $(f_2)$,
\eqref{a2}, \eqref{a22} and \eqref{a3},
we find that
$$|\psi_n(t,\alpha)|\le C_3\Bigl(\frac{Q(t)}{q(t)}\Bigr)^{(2-p)/(p-1)}\frac{Q(t)}{q(t)}
(1+\int_0^t|\psi_n(s,\alpha)|ds),$$ implying that \be\label{a4}
|\psi_n(t,\alpha)|\le C_3t^{p'-1}(1+\int_0^t|\psi_n(s,\alpha)|ds).
\ee
It follows now from Gronwall's lemma that $\psi_n(t,\alpha)$
is uniformly bounded in $[0,t^*]$, and thus from \eqref{a2} and
\eqref{a22} also $\varphi_n(t,\alpha)$, $\frac{\partial}{\partial
t}\varphi_n(t,\alpha)$, and $\Theta_n(t,\alpha)$ are uniformly
bounded in $[0,t^*]$. Also, since
\begin{eqnarray*}
-\frac{\partial}{\partial t}\Theta_n(t,\alpha)=\frac{f(v(t,\alpha+h_n))-f(v(t,\alpha))}
{v(t,\alpha+h_n)-v(t,\alpha)}\varphi_n(t,\alpha)\qquad\qquad\qquad\qquad\qquad\qquad\\
\qquad\qquad\qquad\qquad\qquad\qquad -\frac{q_t(t)}{q^2(t)}
\int_0^tq(s)\frac{f(v(s,\alpha+h_n))-f(v(s,\alpha))}
{v(s,\alpha+h_n)-v(s,\alpha)}\varphi_n(s,\alpha)ds,
\end{eqnarray*}
we find that
$$\Bigm|\frac{\partial}{\partial t}\Theta_n(t,\alpha) \Bigm|\le C_4(1+\frac{tq_t(t)}{q(t)}).$$
From \eqref{a00}, we see that
$$\frac{tq_t(t)}{q(t)}=\frac{1}{p'}\Bigl(\frac{n-1}{p-1}p+\frac{r(t)K'(r(t))}{K(r(t))}\Bigr)
\frac{\int_0^{r(t)}K^{1/p}(\tau)d\tau}{r(t)K^{1/p}(r(t))}.$$ Also,
from $(K_1)$ and the fact that
$(p+rK'(r)/K(r))K^{1/p}=p(rK^{1/p})'$, we obtain that
$$prK^{1/p}(r)\ge \Bigl(p+\frac{rK'(r)}{K(r)}\Bigr)\int_0^rK^{1/p}(\tau)d\tau,$$
and thus
\begin{eqnarray*}
\frac{tq_t(t)}{q(t)}&\le&
\frac{p}{p'}\Bigl(\frac{n-1}{p-1}p+\frac{r(t)K'(r(t))}{K(r(t))}\Bigr)\Bigl(p+\frac{r(t)K'(r(t))}{K(r(t))}\Bigr)^{-1}\\
&\le& (n-p)p\Bigl(p+\frac{r(t)K'(r(t))}{K(r(t))}\Bigr)^{-1}+p-1.
\end{eqnarray*}
Since the last term in this inequality is increasing by $(K_1)$, we conclude that $tq_t/q$ is bounded in $[0,t^*]$.

Hence, from Arzela-Ascoli's Theorem, we conclude that $\varphi_n(\cdot,\alpha)$ and $\Theta_n(\cdot,\alpha)$
converge uniformly in $[0,t^*]$ (up to a subsequence) to continuous functions
$\varphi(\cdot,\alpha)$ and $\Theta(\cdot,\alpha)$
respectively. Moreover, by Lebesgue's dominated convergence theorem, it holds that
$$-\Theta(t,\alpha)=\frac{1}{q(t)}\int_0^tq(s)f'(v(s,\alpha))\varphi(s,\alpha)ds.$$

On the other hand, using that
$$\Theta_n(t,\alpha)=\frac{\phi_p(v_t(t,\alpha+h_n))-\phi_p(v_t(t,\alpha))}
{v_t(t,\alpha+h_n)-v_t(t,\alpha)}\psi_n(t,\alpha),$$
we deduce that for each $t\in(0,t^*]$,
$$\lim_{n\to\infty}\psi_n(t,\alpha)=\frac{\Theta(t,\alpha)}{(p-1)|v_t(t,\alpha)|^{p-2}}:=\psi(t,\alpha),$$
hence again by Lebesgue's dominated convergence theorem we obtain
that \be\label{a5} \varphi(t,\alpha)=1+\int_0^t\psi(s,\alpha)ds.
\ee We will see next that the solution to \be\label{sys}
\begin{gathered}
-\Theta(t,\alpha)=\frac{1}{q(t)}\int_0^tq(s)f'(v(s,\alpha))\varphi(s,\alpha)ds\\
\varphi(t,\alpha)=1+\frac{1}{p-1}\int_0^t\Theta(s,\alpha)|v_t(s,\alpha)|^{2-p}ds,
\end{gathered}
\ee
is unique, and therefore
 the complete sequence $\{\varphi_n\}$ converges to $\varphi$, and the complete
sequences $\{\Theta_n\}$, $\{\psi_n\}$ converge to $\Theta$ and
$\psi$ respectively. Indeed, the only delicate case for the
uniqueness is $p>2$, so let us assume that $p>2$ and that
\eqref{sys} has two solutions $(\varphi_1,\Theta_1)$ and
$(\varphi_2,\Theta_2)$. Substracting the two corresponding
equations we obtain that \be\label{sys2}
\begin{gathered}
-(\Theta_1(t,\alpha)-\Theta_2(t,\alpha))=\frac{1}{q(t)}
\int_0^tq(s)f'(v(s,\alpha))(\varphi_1(s,\alpha)-\varphi_2(s,\alpha))ds\\
\varphi_1(t,\alpha)-\varphi_2(t,\alpha)=
\frac{1}{p-1}\int_0^t(\Theta_1(s,\alpha)-\Theta_2(s,\alpha))|v_t(s,\alpha)|^{2-p}ds.
\end{gathered}
\ee
Using \eqref{a0} to bound $|v_t|$ and the fact that $Q(t)\le tq(t)$ for all $t>0$, we find from the
second equation in \eqref{sys2} that
\be\label{b1}
|\varphi_1(t,\alpha)-\varphi_2(t,\alpha)|\le C_1\int_0^t|\Theta_1(s,\alpha)-\Theta_2(s,\alpha)|s^{(2-p)/(p-1)}ds,
\ee
for some positive constant $C_1$. Also, from the first equation in \eqref{sys2} we see that
$$|\Theta_1(s,\alpha)-\Theta_2(s,\alpha)|\le C_2||\varphi_1(\cdot,\alpha)-\varphi_2(\cdot,\alpha)||_{[0,t]}s,$$
for some positive constant $C_2$, where we have denoted
$$||\varphi_1(\cdot,\alpha)-\varphi_2(\cdot,\alpha)||_{[0,t]}=
\sup\limits_{s\in[0,t]}|\varphi_1(s,\alpha)-\varphi_2(s,\alpha)|,$$
and  therefore, replacing into  \eqref{b1} we obtain \be\label{b2}
|\varphi_1(t,\alpha)-\varphi_2(t,\alpha)|\le
C_3||\varphi_1(\cdot,\alpha)-\varphi_2(\cdot,\alpha)||_{[0,t]}
t^{p'}, \ee for some positive constant $C_3$. We conclude then
that if $t$ is small enough, for example $t\in[0,\bar t]$, with
$C_3\bar t^{p/(p-1)}\le1/2$, then $\varphi_1\equiv\varphi_2$ in
$[0,\bar t]$. The fact that $\varphi_1\equiv\varphi_2$ in
$[0,t^*]$ follows from the standard theory of ordinary
differential equations, so we omit it.

Since
these arguments apply to every sequence $\{h_n\}\to 0$, it follows that $v$ and $\phi_p(v_t)$
are differentiable with respect to $\alpha$ for every $t\in[0,t^*]$
and their derivatives are respectively given
by $\varphi$ and $\Theta$. Also, $v_t$ is differentiable with respect to $\alpha$ for $t\in(0,t^*]$ and
its derivative is given by $\psi$.

The rest of the proof is exactly the same as in \cite{cfe2}: in
order to continuously extend the functions $\varphi$, $\Theta$ and
$\psi$ to $\underline{\mathcal O}$, we replace the  function $f$
in \eqref{pghms-trans2} by $\hat f$ in such a way that $f\equiv
\hat f$ in $[u_0,\infty)$ and $\hat f$ is Lipschitz continuous on
$[0,\infty)$. We can repeat the arguments used above in an
interval $[\delta, T]$ containing $t_0=t(r_0(\alpha))$, with
$\delta>0$, so that we do not need the  estimate \eqref{a1} used
to  bound from below $v_t(t)$ near $t=0$. The solution $\hat v$ of
the new problem satisfies $\hat v\equiv v$ in $[\delta,t_0)$.

The lemma follows by returning to the original variable $r$: $u$
and $K^{-1/p'}\phi_p(u')$ are differentiable with respect to
$\alpha$ in $[0,r_0)$, and these derivatives can be continuously
extended to $[0,r_0]$; $u'(r,\alpha)$ is differentiable with
respect to $\alpha$ in $(0,r_0)$ and
$$\frac{\partial
u'}{\partial\alpha}(r,\alpha)=\psi(r,\alpha)K^{1/p}(r)$$ can be
continuously extended to $[0,r_0]$. Then, from \eqref{a5} we
conclude that
$$\varphi(r,\alpha)=1+\int_0^r\frac{\partial
u'}{\partial\alpha}(\rho,\alpha)d\rho,$$ and hence,
$\phi'(r,\alpha)=(\partial u'/\partial\alpha)(r,\alpha)$ for
$r>0$. Thus we obtain that
$$r^{n-1}|u'(r,\alpha)|^{p-2}\varphi'(r,\alpha)=\frac{\partial
}{\partial\alpha}(r^{n-1}\phi_p(u'(r,\alpha)))$$ and, using
\eqref{sys}, that
$$\frac{\partial
}{\partial\alpha}(r^{n-1}\phi_p(u'(r,\alpha)))=-\int_0^r\rho^{n-1}K(\rho)f'(u(\rho,\alpha))\varphi(\rho,\alpha)d\rho$$
for $r>0$, and thus
$\lim\limits_{r\to0^+}r^{n-1}|u'(r)|^{p-2}\varphi'(r,\alpha)=0$
and $\varphi$ satisfies \eqref{varphi-eq} in $(0,r_0]$.

Finally, from \eqref{sys},  it can be seen that $\varphi$
satisfies
$$\varphi(r,\alpha)-1=\frac{1}{p-1}\int_0^ru'(\rho,\alpha)
\left(\frac{\int_0^\rho\tau^{n-1}K(\tau)f'(u(\tau,\alpha))\varphi(\tau,\alpha)d\tau}{\rho^{n-1}|u'(\rho,\alpha)|^{p-1}}
\right)d\rho,$$ and thus, after some computations involving L'H\^
opital's rule, we obtain that $\varphi$ is continuously
differentiable with respect to $r$ at $r=0$, and
$$\varphi'(0)=0=\frac{\partial u'}{\partial\alpha}(0,\alpha).$$

\end{document}